\documentclass[twoside, 11pt]{article}
\usepackage{mathrsfs}

\usepackage{mathrsfs,amsfonts,amsmath}
\RequirePackage{ifthen,calc}
\usepackage{theorem}
\usepackage[dvips]{color}

 \catcode`@=11
 \def\@evenhead{\hbox to\textwidth{\footnotesize\rm\thepage \hfill
  {\it }}} 

 \def\@oddhead{\hbox to \textwidth{\footnotesize{\it
  Small deviation for Random walk with random environment in time } \hfill\thepage}}

 \renewcommand{\section}{\makeatletter
 \renewcommand{\@seccntformat}[1]{{\csname the##1\endcsname.}\hspace{0.45em}}
 \makeatother \@startsection
{section}
{1}
{0pt}
{\baselineskip}
{0.5\baselineskip}
{\normalsize\bfseries\mathversion{bold}}}

\catcode`@=12
\newcommand\ack{\section*{Acknowledgement}}

\newtheorem{thm}{\noindent Theorem}[section]
\newtheorem{lem}{\noindent Lemma}[section]
\newtheorem{cor}{\noindent Corollary}[section]

\theoremheaderfont{\normalfont\bfseries} \theorembodyfont{\slshape}
{\theorembodyfont{\rmfamily}}
{\theorembodyfont{\rmfamily}}
{\theorembodyfont{\rmfamily}}
{\theorembodyfont{\rmfamily}}
{\theorembodyfont{\rmfamily}}
{\theorembodyfont{\rmfamily}\newtheorem{rem}{\noindent Remark}[section]}
{\theorembodyfont{\rmfamily}}
{\theorembodyfont{\rmfamily}}
{\theorembodyfont{\rmfamily}}

 \def\beqlb{\begin{eqnarray}}\def\eeqlb{\end{eqnarray}}
 \def\beqnn{\begin{eqnarray*}}\def\eeqnn{\end{eqnarray*}}
 \numberwithin{equation}{section}
\def\qed{\hfill$\square$\smallskip}

\def\bfE{{\mathbb{E}}}
\def\bfP{{\mathbb{P}}}
\def\bfR{{\mathbb{R}}}

\def\bfN{{\mathbb{N}}}

\begin{document}
\title{\bf Small deviation for Random walk with random environment in time
}
\author{ You Lv\thanks{Email: youlv@mail.bnu.edu.cn }
\\ \small School of science of mathematics, Beijing Normal University,
\\ \small Beijing 100875, P. R. China.
}
\date{}

\maketitle


\noindent\textbf{Abstract}: We give the random environment version of Mogul'ski\v{\i} estimation in quenched sense.~
 Assume that $\{\mu\}_{n\in\bfN}$ (called environment) is a sequence of i.i.d. random probability measures on $\bfR.$~ Let $\{X_n\}_{n\in\bfN}$ be a sequence of independent random variables, where $X_n$ has law $\mu_n.$ We set $S_n=\sum_{i=1}^{n}X_i.$ Under some integrability conditions, we show that on the log scale, for any power function $f$. the decay rate of $$\bfP_\mu(\forall_{0\leq i\leq n} S_{f(n)+i}\in[g(i/n)n^{\alpha},h(i/n)n^{\alpha}]|S_{f(n)}=x)$$ is $e^{-cn^{1-2\alpha}}$ almost surely as $n\rightarrow+\infty$, where $c>0,\alpha\in(0,\frac{1}{2}),$ $g,h\in\mathcal{C}[0,1]$ (the set of all continuous functions defined on $[0,1]),$ $g(s)<h(s), \forall s\in[0,1],$ and $x\in(g(0),h(0)).$ The main result of this paper is also a basic tool in the researching of Branching random walk in random environment with selection.

\smallskip


\smallskip

\noindent\textbf{2000 Mathematics Subjects
Classification:} 60G50

\section{Introduction}
We introduce a model called random walk with random environment in time. Let $\mu:=\{\mu_n\}_{n\in\bfN}$ be a sequence of i.i.d. random probability measures on $\bfR.$  Conditionally on a realization of $\mu,$ we introduce a sequence of independent random variables  $\{X_n\}_{n\in\bfN}$ such that the law of $X_n$ is $\mu_n$ for every $n\in\bfN.$  We denote by $S_n:=\sum_{i=1}^{n}X_i$ the random walk with random environment in time. We write $\bfP_{\mu}$ for the law of  $\{S_n\}_{n\in\bfN}$ conditionally on the sequence $\{\mu_n\}_{n\in\bfN},$ $\bfP$ for the joint law of $\{S_n\}_{n\in\bfN}$ and the environment $\mu.$ The corresponding expectations are denoted by $\bfE_{\mu}$ and $\bfE$ respectively. As the conventional notations, we denote $\bfP^x_{\mu}(\cdot):=\bfP^x_{\mu}(\cdot|S_0=x)$ and $\bfE^x_{\mu}(\cdot):=\bfE^x_{\mu}(\cdot|S_0=x).$

 This paper focuses on the small deviation probabilities
of $S_n.$  More precisely, if we let \beqlb\label{sec-0.0}H_{n,f}:=\big\{\forall_{0\leq i\leq n} S_{f(n)+i}\in[g(i/n)n^{\alpha},h(i/n)n^{\alpha}]\big\},\eeqlb we want to obtain the asymptotical behavior of $\bfP_{\mu}(H_{n,f}|S_{f(n)}=x),$ where $g(s), h(s)$ are two continuous functions defined on $[0,1]$ such that $g(s)<h(s)$ for any $s\in[0,1],$ $x\in (g(0), h(0)), \alpha\in(0,\frac{1}{2})$ and $f(n)$ meets the assumption (H4) which is listed in section 2. The starting time $f(n)$ should be emphasized because~$\{S_n\}_{n\in\bfN}$ is a time-inhomogeneous random walk conditionally on a realization of $\mu$.  If the environment space is degenerate, it is the so-called Mogul'ski\v{\i}
estimation. In 1974, Mogul'ski\v{\i} \cite{Mog1974} first showed the order decay of $\bfP^x(H_{n,0})$ is $e^{-cn^{1-2\alpha}}$ when $S_n$ is an i.i.d. centered random walk with finite variance. As a basic tool, it has attracted
many scholars' attention because of its wide application in probability and statistics models. Gantert, Hu and
Shi \cite{GHS2011} gave the triangular version of  Mogul'ski\v{\i} estimation. Namely, let $S^{n}_i=\sum_{k=1}^{i}X^{n}_k, 1\leq i\leq n,$ where $\{X^{n}_i\}_{i\in\{1,2,...n\}}$ is an i.i.d. random sequence for any fixed $n$ and $X^{n}_1$ meets some integrability conditions. They discussed the probability $\bfP^x(H_{n,0})$ if we replace the $S_i$ by $S^n_i$ in the definition (1.1).
For the sum of independent but not necessarily identically distributed random variables, the Mogul'ski\v{\i} estimation has been investigated in \cite{Mal2015} . Mallein utilized the
Mogul'ski\v{\i} estimation of time-inhomogeneous version to study the Maximal displacement of a branching random walk with time-inhomogeneous environment. Shao \cite{QM1995} discussed more generalized small deviation probabilities for the sum of independent centred random variables.  For more
results on small deviation probabilities, see \cite{L2006}.

For a random walk $S_n$ with random environment in time, \cite{MM2015} and \cite{MM2016} have obtained some properties on its first hitting time. For example, the Random Ballot Theorem in \cite{MM2015} showed the decay rate of $\bfP^0_{\mu}(\forall_{i\leq n} S_i\geq -1)$  is $n^{-\gamma}$ almost surely as $n\rightarrow +\infty$, where $\mu$ is the i.i.d. random environment and $\gamma$ is a positive constant depending on $S_n.$  $\cite{MM2016}$ can be seen as an application of the Random Ballot Theorem in branching random walk with
random environment in time. The main goal of this paper is to study the limit behavior of probability $\bfP^x_\mu(H_{n,f})$ in quenched sense, which can be seen
as an application of \cite{LY2018} in random walk with random environment in time. Moreover, the main result in this paper is also a essential tool when we study the Branching random walk in random environment with selection, see \cite{LY201809}.

The rest of this paper is organized as follows. We state the main theorem and corollary in section 2. The proof is given in section 3.

\section{Main result}
 First, we need some basic assumptions for $S_n$ throughout this paper. Denote $$M_{n}:=\bfE_\mu(S_n), ~U_{n}:=S_n-\bfE_\mu(S_n),~\Gamma_n:=\bfE_\mu(U^2_n)=\bfE_\mu(S^{2}_n)-M^2_{n}.$$
\begin{itemize}
  \item[(H1)] $\bfE M_1=0,\sigma^2_A:=\bfE(M_1^2)\in[0,+\infty),\sigma^2_Q:=\bfE(U_1^2)=\bfE(\Gamma_1)\in(0,+\infty).$
  \item[(H2)] There exists $\lambda_1>0,$ such that $\bfE(e^{\lambda_1 |M_{1}}|)<+\infty.$
  \item[(H3)] There exist $\lambda_2, \lambda_3>0$ such that $\bfE_\mu(e^{\lambda_2|U_1|})\leq \lambda_3$ almost surely.
  \item[(H4)] $\{r_n\}_{n\in\bfN}$ is a positive sequence such that for any $\varrho>0, \lim\limits_{n\rightarrow+\infty}\frac{n^{^\varrho}}{r_{_n}}=0.$ $f(n)$ is an positive integer-valued function such that for any $\kappa>0, \lim\limits_{n\rightarrow+\infty}\frac{f(n)}{e^{n^{\kappa}}}=0.$
\end{itemize}
\begin{rem}
(H2) and (H3) are the technical assumptions for using KMT theorem \cite[Corollary 2.2]{Lif2000} and Sakhanenko theorem \cite[Theorem 3.1]{Lif2000}. According to (H3), we deduce that for all $\lambda'_2 <\lambda_2,$  there exists $C>0$ such that for all $x\geq0,$ $x^3e^{\lambda'_2x} \leq Ce^{\lambda_2x}.$ That implies the Sakhanenko parameter \cite[Page 7]{Lif2000} exists. Hence the conditions of Sakhanenko theorem are met.
\end{rem}

Now we introduce the main result of this paper. In order to make our theorem more flexible and have wider application, we add some other events to the lower bound. Let $\xi_i$ be a positive random variable whose law is only determined by the $i-$th element $\mu_i$ in a realization of environment $\mu.$ Therefore,
conditionally on a given environment realization $\mu,$  $\{\xi_i\}_{i\in\bfN}$ is an independent positive random sequence since $\mu$ is i.i.d.. Moreover, for any measurable function $\eta$, $\{\bfE_{\mu}(\eta(\xi_i))\}_{i\in\bfN}$ is an i.i.d. random sequence in the environment space. We also need to introduce a function $\gamma:\bfR \rightarrow \bfR_+$ given by the following limit
$$\lim \limits_{t\rightarrow +\infty}\frac{-\ln \bfP(\forall_{s\leq t} |B_s-\beta W_s|\leq 1/2|W)}{t}=\gamma(\beta),~~~\rm{a.s.},$$
where $B,W$ are independent standard Brownian motions. This function are studied in detail in \cite{LY2018}.
\begin{thm}
Let $\alpha\in(0,\frac{1}{2}), a<a_0\leq b_0<b,a\leq a'<b'\leq b, f'(n):=f(n)+n.$ Under the assumption $\bfE(\xi_1)<+\infty$ and (H1)-(H4), we have, almost surely
$$\limsup \limits_{n\rightarrow +\infty}~\sup\limits_{x\in\bfR}\frac{\ln \bfP_\mu
(\forall_{f(n)\leq i\leq f'(n)} an^\alpha\leq S_{i}\leq bn^\alpha|S_{f(n)}=x)}{n^{1-2\alpha}}\leq -\frac{\sigma^2_Q}{(b-a)^2}\gamma(\frac{\sigma_A}{\sigma_Q}),$$
\beqnn&&\liminf \limits_{n\rightarrow +\infty}\inf\limits_{x\in[a_0 n^{\alpha}, b_0 n^{\alpha}]}\frac{\ln \bfP_\mu
\left(\begin{split}\forall_{f(n)\leq i\leq f'(n)} S_{i}\in [a n^\alpha,b n^\alpha],\\S_{f'(n)}\in [a'n^\alpha,b'n^\alpha],\xi_i\leq r_n\end{split}\Bigg|S_{f(n)}=x\right)}{n^{1-2\alpha}}
\\&&~~~~~~~~~~~~~~~~~~~~~~~~~~~~~~~~~~~~~~~~~~~~~~~~~~~~~~~~~~~~~~~~~~~~~~~~\geq -\frac{\sigma^2_Q}{(b-a)^2}\gamma(\frac{\sigma_A}{\sigma_Q}).\eeqnn
\end{thm}
\begin{cor}
Under the assumption $\bfE(\xi_1)<+\infty$ and (H1)-(H4), let $g(s), h(s)$ be two continuous functions on $[0,1]$ and $g(s)<h(s)$ for any $s\in [0,1].$ We set $g(0)< a_0\leq b_0 <h(0), g(1)\leq a'<b'\leq h(1)$
~and~ $C_{g,h}:=\int_{0}^{1}\frac{1}{[h(s)-g(s)]^2}ds.$ Then for any $\alpha\in (0,\frac{1}{2}),$ we have, almost surely,
$$\limsup \limits_{n\rightarrow +\infty}\sup\limits_{x\in\bfR}\frac{\ln\bfP_\mu
\Big(\forall_{0\leq i\leq n}\frac{S_{f(n)+i}}{n^\alpha}\in \big[g\big(\frac{i}{n}\big),h\big(\frac{i}{n}\big)\big]\Big|S_{f(n)}=x\Big)}{n^{1-2\alpha}}\leq -C_{g,h}\sigma^2_Q\gamma(\frac{\sigma_A}{\sigma_Q}),$$
\beqnn&&\liminf \limits_{n\rightarrow +\infty}\inf\limits_{x\in[a_0 n^{\alpha}, b_0 n^{\alpha}]}\frac{\ln \bfP_\mu
\Bigg(\begin{split}\forall_{0\leq i\leq n} S_{f(n)+i}\in [g(i/n)n^\alpha,h(i/n)n^\alpha],\\ S_{f'(n)}\in [a'n^\alpha,b'n^\alpha],~\xi_{i+f(n)}\leq r_n \end{split} \Bigg|S_{f(n)}=x\Bigg)}{n^{1-2\alpha}}
\\&&~~~~~~~~~~~~~~~~~~~~~~~~~~~~~~~~~~~~~~~~~~~~~~~~~~~~~~~~~~~~~~~~~~~~~~\geq -C_{g,h}\sigma^2_Q\gamma(\frac{\sigma_A}{\sigma_Q}).
\eeqnn
 \end{cor}

 \begin{rem}
The Random ballot theorem in \cite{MM2015} tells us that the decay rate of $\bfP_\mu^x(\forall_{i\leq n}S_i\geq -1)$ is only depend on $\frac{\sigma_A}{\sigma_Q}.$ But here, Theorem 2.1 tells us that the decay rate of small deviation probability will depend on $\sigma_A$ and $\sigma_Q.$ We can intuitively comprehend this phenomenon by recalling the following two basic conclusions. For a standard Brownian motion $(Z_t, t\geq 0)$ with parameters $\bfE(Z_t)=0$ and  $\bfE(Z^2_t)=\sigma^2t,$ it is well known that $$\lim\limits_{t\rightarrow+\infty}\frac{\ln \bfP(\forall_{s\leq t}Z_s\geq -1)}{\ln t}=-\frac{1}{2},~\lim\limits_{t\rightarrow+\infty}\frac{\ln \bfP(\forall_{s\leq t}|Z_s|\leq 1)}{t}=-\frac{\pi^2\sigma^2}{8}.$$   In our model, $\sigma^2_Q$ plays the role as the variance $\sigma^2$ of $Z$ in some extent.

According to \cite[Theorem 2.1]{LY2018},  we know $\gamma(x)$ is a convex and strictly increasing on $[0,+\infty).$ It implies that compared with a time homogeneous random walk with expectation $0$ and variance $\sigma_Q,$ the decay rate of the small deviation probabilities in our model is larger as soon as almost surely $\sigma_A>0$. Moreover, Corollary 2.1 is agree with the Mogul'ski\v{\i} estimation in \cite{Mog1974} when the environment is degenerate since $\gamma(0)=\frac{\pi^2}{2}.$

\end{rem}
\section{Proof of main result}
The main idea of the proof is coupling $S_n$ with two Brownian motions by KMT theorem \cite[Corollary 2.2]{Lif2000} and Sakhanenko theorem
\cite[Theorem 3.1]{Lif2000}. These coupling method has been used in \cite{MM2015}. To using that, we should first decompose $S_{n}$ as the following way.
Noting that $S_{n}$ is the sum of $M_{n}$ and $U_{n}$, ${M_{n}}$ is an i.i.d. centred random walk since the random environment is i.i.d.. When we give a realization of $\{\mu_{n}\}_{n\in\bfN}$, $U_{n}=\sum\limits_{i=0}\limits^{n}\left(X_i-E_{\mu}(X_i)\right)$ is the sum of centred independent (but not necessarily identically distributed) random variables.
$\{\Gamma_n\}_{n\in\bfN}$ is sums of i.i.d. non-negative random variables and $\bfE(\Gamma_1)=\sigma^2_Q.$
  ~The concrete detail of how to implement the coupling is totally same as \cite[Page30-31]{MM2015}, so here we only
 give the final form after coupling. According to the KMT theorem, we can construct a standard Brownian
motion $W$
to satisfy that there exists positive constants $\lambda_0,C_0,D$ such that
\beqlb\label{sec-1}\forall x\geq 0,~ \forall n\in \bfN,~ \bfP(\Delta'_{n}\geq D\ln n+x)\leq C_0e^{-\lambda_0 x},\eeqlb where $\Delta'_{n}:=\sup\limits_{i\leq n}|M_i-\sigma_A W_{i}|.$
Next, by Sakhanenko Theorem we approximate $\{U_n\}$ by a Brownian motion $B$ which is independent of $W$. If we set $\Delta_{n}:=\sup\limits_{i\leq n}|U_i-B_{\Gamma_i}|,$ then there exist positive
constants $C$ and $\lambda$ such that
$$\bfE_{\mu}\left[e^{C\lambda\Delta_{n}}\right]\leq 1+\lambda\sqrt{\Gamma_n}~,\rm{~~a.s.}.$$

\begin{lem}
Let $B, W$ be two independent standard Brownian motions, $W_0=0$. We denote $\bfP^x(\cdot|W):=\bfP(\cdot|W,B_0=x),$ which is the conditional probability for $B$ when we give a realization of $W$ and the starting point of $B.$ Assume that $a<a_0\leq b_0<b, a\leq a'<b'\leq b$. For any $t,p,m>0$ and positive function $f,$
denote
\beqnn \overline{X}_{n,t}:=-\ln\inf\limits_{x\in [a_0,b_0]}\bfP^x\left(\begin{split}\forall_{s\leq t} ~B_{s}+\beta n^{-\alpha}\left(W_{f(n)+sn^{2\alpha}}-W_{f(n)}\right)\in[a, b],
\\ ~~~~B_{t}+\beta n^{-\alpha}\left(W_{f(n)+tn^{2\alpha}}-W_{f(n)}\right)\in\left[a',b'\right]\end{split}\Bigg|W\right).\eeqnn Then we can find a $q>2, N(q)\in \bfN$ such that for every $m\geq N(q),$ we have \beqlb\label{sec-2}\bfP\left(\overline{X}_{n,t}\geq e^{m^p}\right)\leq m^{-q}.\eeqlb
\end{lem}
\noindent{\bf Proof of lemma 3.1}
Notice that for any $n,$  $X_{n,t}$ has the same distribution as
$$-\ln\inf\limits_{x\in [a_0,b_0]}\bfP^x\big(\forall_{s\leq t} B_{s}+\beta W_{s}\in[a, b],
B_{t}+\beta W_{t}\in[a',b']|W\big).$$
Therefore, from \cite[Theorem 3.1]{LY2018} we can know that (3.2) is true.   \qed

Now, let us show the lower bound first.

\noindent{\bf Proof of Theorem 2.1: the lower bound}

First we define some important events. Let $D\in\bfR^+, ~T=\lfloor Dn^{2\alpha}\rfloor,~ K=\left\lfloor \frac{n}{T} \right\rfloor,$ $f_k(n):=f(n)+kT.$
Recalling that we use Brownian motion $\{W_s,s\geq 0\}$ to approximate $\{M_i, i\geq 0\}$.  Using Brownian motion $\{B_s, s\geq 0\}$ to approximate
$\{U_i, i\geq 0\}.$ Choose $\rho,\varsigma,u>0$ such that $0<\frac{\alpha+\varsigma}{2}<\rho<\alpha, u\in(0,1).$ We define
$$A_{n}:=D_n\cap E_n\cap F_n\cap G_{n}\cap I_n,$$
Where $$D_n:=\{\Delta'_{f(n)+n}\leq n^\rho\},~~ \Delta'_{f(n)+n}:=\sup\limits_{i\leq f(n)+n}\left|M_i-\sigma_A W_{i}\right|;$$
$$E_n:=\left\{\sup_{s\in[i,i+1)} \sigma_{A}|W_s-W_i|\leq n^{\rho}, \forall f(n)\leq i\leq f(n)+n\right\};$$
$$F_n:=\cap_{k=0}^{K}\big\{\forall_{i\leq T}|\Gamma_{f_k(n)+i}-\Gamma_{f_k(n)}-i\sigma^2_Q|\leq n^{\alpha+\varsigma}\big\},~I_n:=\left\{{\frac{\sum_{i=1}^{f(n)+n}\bfE_\mu(\xi_i)}{\bfE(\xi_1)(f(n)+n)}\leq 2}\right\}.$$
We will give the definition of $G_{n}$ postponed.

By (3.1) and the Borel-Cantelli 0-1 law, we have $\mathbf{1}_{D_n}\rightarrow 1,$ a.s. as $n\rightarrow+\infty.$

By basic calculation of Brownian motion and the Borel-Cantelli 0-1 law we can also know $\mathbf{1}_{E_n}\rightarrow 1,$ a.s. as $n\rightarrow+\infty.$

Since $\{\bfE_\mu(\xi_i)\}_{i\in\bfN}$ is a sequence of i.i.d. random variables with finite mean $\bfE(\xi_1),$ according to the Borel-Cantelli 0-1 law and the large deviation principle, it is true that $\mathbf{1}_{I_n}\rightarrow 1,$ a.s. as $n\rightarrow+\infty.$

Choose a positive even number $p$ such that $p>\frac{3-2\alpha}{\varsigma}.$ We can see $\{|\Gamma_i-i\sigma^2_Q|, i\geq 0\}$ is a submartingale taking
non-negative value. So using the Doob's martingale inequality we have
\beqnn(K+1)\bfP\Big(\sup_{0\leq i \leq T}|\Gamma_i-i\sigma^2_Q|\geq n^{\alpha+\varsigma}\Big)&\leq& (K+1)\frac{\bfE(\sup_{0\leq i\leq T}|\Gamma_i-i\sigma^2_Q|^p)}{n^{p(\alpha+\varsigma)}}
\\&\leq& (K+1)\Big(\frac{p}{p-1}\Big)^p\frac{\bfE(\Gamma_T-T\sigma^2_Q)^p}{n^{p(\alpha+\varsigma)}}
\\&=& (K+1)\Big(\frac{p}{p-1}\Big)^p\frac{C_1\mathcal{C}_{T}^{\frac{p}{2}}+o(n^{2\alpha\cdot\frac{p}{2}})}{n^{p(\alpha+\varsigma)}},  \eeqnn
where $\mathcal{C}$ is the combinatorial number, $C_1=\big(\bfE(|\Gamma_1-\sigma_Q|^2)\big)^{\frac{p}{2}}.$ From (H3), we know $C_1<+\infty$, which implies that
 $$(K+1)\bfP\left(\sup_{0\leq i \leq T}|\Gamma_i-i\sigma^2_Q|\geq n^{\alpha+\varsigma}\right)<\frac{K}{n^{3-2\alpha}}<\frac{1}{n^2}.$$Hence $\mathbf{1}_{F_n}\rightarrow 1,$ a.s. as $n\rightarrow+\infty$ by the Borel-Cantelli 0-1 law.

The event $A_n$ is focus on process $M$ and Brownian motion $W.$ The following events $C_{n,k}$ and $D_n$ are related with process $U$ and Brownian motion $B.$
Recalling that $f_k(n):=f(n)+kT$. Let
$$\Delta_{n,k}:=\sup\limits_{i\leq T}\left|U_{f_k(n)+i}-U_{f_k(n)}-(B_{\Gamma_{f_k(n)+i}}-B_{\Gamma_{f_k(n)}})\right|~,~C_{n,k}:=\{\Delta_{n,k}\leq n^{\rho}\}.$$
By Sakhanenko theorem, almost surely, there exist $C, \lambda>0$ such that
\beqlb\label{sec-3}\bfP_{\mu}(C^c_{n,k})&\leq& \bfP_{\mu}\left(\sup\limits_{i\leq T}|U_{f_k(n)+i}-B_{\Gamma_{f_k(n)+i}}|+|U_{f_k(n)}-B_{\Gamma_{f_k(n)}}|\geq n^\rho\right)\nonumber
\\&\leq& \bfP_{\mu}\left(\sup\limits_{i\leq T}|U_{f_k(n)+i}-B_{\Gamma_{f_k(n)+i}}|\geq \frac{n^\rho}{2}\right)+\bfP_{\mu}\left(|U_{f_k(n)}-B_{\Gamma_{f_k(n)}}|\geq \frac{n^\rho}{2}\right)\nonumber
\\&\leq& 2\bfP_{\mu}\left(\Delta_{f(n)+n}\geq \frac{n^{\rho}}{2}\right)\leq \frac{2\bfE_{\mu}\big(e^{C\lambda\Delta_{f(n)+n}}\big)}{e^{\frac{C\lambda n^{\rho}}{2}}}
\leq \frac{2+2\lambda\sqrt{\Gamma_{f(n)+n}}}{e^{\frac{C\lambda n^{\rho}}{2}}}.
\eeqlb
Denote
$$D_{n}:=\left\{\forall i\leq T,\forall |t|\leq n^{\alpha+\varsigma}+\sigma^2_Q, |B_{i\sigma^2_{Q}}-B_{i\sigma^2_{Q}+t}|\leq \frac{1}{2}n^\rho\right\}.$$
According to the basic property of Brownian motion, we have
\beqlb\label{sec-4}\bfP(D^c_n)&\leq&2T\bfP\left(\sup_{t\in[0,n^{\alpha+\varsigma}+\sigma^2_Q]}|B_t|>n^{\rho}\right)\leq 2T\cdot c_1\exp\left\{\frac{-c_2n^{2\rho}}{n^{\alpha+\varsigma}+\sigma^2_Q}\right\},\eeqlb where $c_1,c_2$ are two positive constants.
Notice that on the event $F_n\cap D_{n},$ for each $s\in[i,i+1),$ we have
$$\left|B_{s\sigma^2_Q}-B_{\Gamma_{f_k(n)+i}-\Gamma_{f_k(n)}}\right|\leq \left|B_{s\sigma^2_Q}-B_{i\sigma^2_Q}\right|+\left|B_{i\sigma^2_Q}-B_{\Gamma_{f_k(n)+i}-\Gamma_{f_k(n)}}\right|\leq n^{\rho}.$$
Furthermore,we can see $\Gamma_{f(n)+n}\leq K(T\sigma^2_Q+ n^{\alpha+\varsigma})$ on the event $F_n$. Therefore, combining with (3.3), (3.4) and the definition of $I_n,$ for large enough $n$, on the event $A_n,$ we can find a $\upsilon>0$ such that
\beqlb\label{sec-5}\bfP_{\mu}(C_{n,k})+\bfP_{\mu}(D_{n,k})+\frac{\sum_{~~i=1}^{f(n)+n}\bfE_\mu(\xi_i)}{r_n}\leq e^{-n^{\upsilon}}.\eeqlb

Now we begin to do the estimation. Recalling that $f'(n):=f(n)+n$. In this section we assume $a\leq a'<a''<a'''<b'''<b''<b'\leq b, a<a_0\leq b_0<b.$ Let $$\varepsilon_0:=\min\{a'-a, a''-a',a'''-a'',b''-b''',b'-b'',b-b'\},$$ Denote $\sum_{i=f(n)+1}^{f(n)+n}\bfE_{\mu}(\xi_i):=L_n.$ Choose an $\varepsilon\in(0,\varepsilon_0)$ arbitrarily, the following series of inequalities will hold when $n$ is large enough (at least $n^{\alpha-\rho}>\frac{5}{\varepsilon}$).

\beqnn
&&\inf\limits_{x\in[a_0 n^{\alpha}, b_0 n^{\alpha}]}\bfP_\mu
\Big(\forall_{f(n)\leq i\leq f'(n)} n^{-\alpha}S_{i}\in[a,b],n^{-\alpha}S_{f'(n)}\in[a', b'],\xi_i\leq r_n\Big|S_{f(n)}=x\Big)
\\&\geq&\prod\limits_{k=0}\limits^{K-1}\inf\limits_{x\in [a_0 n^\alpha,b_0 n^\alpha]} \bfP_\mu\left(\begin{split}\forall_{i\leq T}~ S_{f_k(n)+i}\in[a n^\alpha, bn^\alpha],\xi_{f_k(n)+i}\leq r_n,\\ S_{f_k(n)+T}\in[a''n^\alpha,b''n^\alpha] ~~~~\end{split}\Bigg|S_{f_k(n)}=x\right)
\\&\times&\inf\limits_{x\in [a'' n^\alpha,b'' n^\alpha]} \bfP_\mu\left(\begin{split}\forall_{i\leq f'(n)-f_K(n)} S_{f_K(n)+i}\in[an^\alpha, b n^\alpha],\\ S_{f'(n)}\in[a'n^\alpha,b'n^\alpha],\xi_{f_K(n)+i}\leq r_n \end{split} \Bigg|S_{f_K(n)}=x\right)
\\&\geq&\prod\limits_{k=0}\limits^{K-1}\left[\inf\limits_{x\in [a_0 n^\alpha,b_0 n^\alpha]} \bfP_\mu\left(\begin{split}\forall_{i\leq T} S_{f_k(n)+i}\in[an^\alpha, bn^\alpha],\\ S_{f_k(n)+T}\in[a''n^\alpha,b''n^\alpha] \end{split} \Bigg|S_{f_k(n)}=x\right)-\frac{L_n}{r_n}\right]
\\&\times&\left[\inf\limits_{x\in [a'' n^\alpha,b'' n^\alpha]} \bfP_\mu\left(\begin{split}\forall_{i\leq f'(n)-f_K(n)}~ S_{f_K(n)+i}\in[a n^\alpha, b n^\alpha],\\ S_{f'(n)}\in[a'n^\alpha,b'n^\alpha]~~\end{split}\Bigg|S_{f_K(n)}=x\right)
-\frac{L_n}{r_n}\right].
\eeqnn
 Let
$$U^{(n,k)}_{s}:=U_{f_k(n)+s}-U_{f_k(n)},~~~ W^{(n,k)}_{s}:=W_{f_k(n)+s}-W_{f_k(n)},$$ $$\Gamma^{(n,k)}_{s}:=\Gamma_{f_k(n)+s}-\Gamma_{f_k(n)},~~~ s\leq T, k\leq K~~~ {\rm and} ~~y:=x n^{-\alpha}.$$ We denote $\bfP_\mu^x(\cdot|W):=\bfP^x(\cdot|W,B_0=x,\mu),
\bfP_\mu(\cdot|W):=\bfP^x(\cdot|W,B_0=0,\mu).$  On the event $A_n,$ for every $k\in [0,K-1]\cap\bfN,$ we have
\beqnn
&& \bfP_\mu\Big(\forall_{i\leq T} n^{-\alpha} S_{f_k(n)+i}\in[a, b],n^{-\alpha}S_{f_k(n)+T}\in[a'',b'']|S_{f_k(n)}=x\Big)
\\&=& \bfP_\mu\left(\begin{aligned}\forall_{i\leq T}& ~~x+U_{f_k(n)+i}-U_{f_k(n)}+M_{f_k(n)+i}-M_{f_k(n)}\in[an^\alpha, bn^\alpha],\\
& x+U_{f_k(n)+T}-U_{f_k(n)}+M_{f_k(n)+T}-M_{f_k(n)}\in[a''n^\alpha,b''n^\alpha]\end{aligned}\right)
\\&\geq& \bfP_\mu\left(\begin{aligned}\forall_{i\leq T} &~ x+U^{(n,k)}_{i}+\sigma_{A}(W_{f_k(n)+i}-W_{f_k(n)})\in[an^\alpha+2n^\rho, bn^\alpha-2n^\rho],
\\ & x+U^{(n,k)}_{T}+\sigma_{A}(W_{f_k(n)+T}-W_{f_k(n)})\in[a'' n^\alpha-2n^\rho,b'' n^\alpha+2n^\rho]\end{aligned}\Bigg|W\right)
\\&\geq& \bfP_\mu\left(\begin{aligned} & \forall_{i\leq T}~ x+B_{\Gamma_{f_k(n)+i}}-B_{\Gamma_{f_k(n)}}+\sigma_{A}W^{(n,k)}_i\in[an^\alpha+3n^\rho,bn^\alpha-3n^\rho],
\\& x+B_{\Gamma_{f_k(n)+T}}-B_{\Gamma_{f_k(n)}}+\sigma_{A}W^{(n,k)}_T\in[a'' n^\alpha+3n^\rho,b'' n^\alpha-3n^\rho],C_{n,k}\end{aligned}\Bigg|W\right)
\\&\geq& \bfP_{\mu}\left(\begin{aligned}\forall_{i\leq T} ~~x+B_{\Gamma^{(n,k)}_{i}}+\sigma_{A}W^{(n,k)}_i\in[an^\alpha+3n^\rho,bn^\alpha-3n^\rho],
\\ x+B_{\Gamma^{(n,k)}_{T}}+\sigma_{A}W^{(n,k)}_T\in[a'' n^\alpha+3n^\rho,b'' n^\alpha-3n^\rho],D_n\end{aligned}\Bigg|W\right)-\bfP_\mu(C^c_{n,k})
\\&\geq& \bfP^x_{\mu}\left(\begin{aligned}\forall_{s\leq T} B_{s\sigma^2_Q}+\sigma_{A}W^{(n,k)}_s\in[an^\alpha+5n^\rho,bn^\alpha-5n^\rho],
\\B_{T\sigma^2_Q}+\sigma_{A}W^{(n,k)}_T\in[a'' n^\alpha+5n^\rho,b'' n^\alpha-5n^\rho]\end{aligned}\Bigg|W\right)-\bfP_\mu(C^c_{n,k})-\bfP_{\mu}(D^c_n)
\\&\geq& \bfP^y\left(\begin{aligned}\forall_{s\leq \frac{T}{n^{2\alpha}}} \sigma_Q B_{s}+\sigma_{A}n^{-\alpha}W^{(n,k)}_{sn^{2\alpha}}\in[a+\varepsilon, b-\varepsilon],
\\ \sigma_Q B_{\frac{T}{n^{2\alpha}}}+\sigma_{A}n^{-\alpha}W^{(n,k)}_{T}\in[a''+\varepsilon,b''-\varepsilon] \end{aligned} \Bigg|W\right)-\bfP_\mu(C^c_{n,k})-\bfP_{\mu}(D^c_n).
\eeqnn
Therefore, combining with (3.5) we know on the event $A_n$, for any $k\in [0,K-1]\cap\bfN,$
\beqnn
&&\inf\limits_{x\in [a_0 n^\alpha,b_0 n^\alpha]} \bfP_\mu\left(\forall_{i\leq T} \frac{S_{f_k(n)+i}}{n^\alpha}\in[a, b],\frac{S_{f_k(n)+T}}{n^\alpha}\in[a',b'],\xi_i\leq r_n|S_{f_k(n)}=x\right)
\\&\geq& 1_{A_n}\inf\limits_{x\in [a_0,b_0]}\bfP^x\left(\begin{aligned}\forall_{s\leq \frac{T}{n^{2\alpha}}} \sigma_Q B_{s}+\sigma_{A}n^{-\alpha}W^{(n,k)}_{sn^{2\alpha}}\in[a+\varepsilon, b-\varepsilon],
\\ \sigma_Q B_{\frac{T}{n^{2\alpha}}}+\sigma_{A}n^{-\alpha}W^{(n,k)}_{T}\in[a''+\varepsilon,b''-\varepsilon] \end{aligned} \Bigg|W\right)
\\&&~~~~~~~~~~~~~~~~~~~~~~~~~~~~~~~~~~~~~-1_{A_n}\left(\bfP_\mu(C^c_{n,k})+\bfP_{\mu}(D^c_n)+\frac{\sum_{i=1}^{f(n)+n}\bfE_\mu(\xi_i)}{r_n}\right)
\\&\geq& 1_{A_n}\left[\inf\limits_{x\in [a_0,b_0]}\bfP^x\left(\begin{aligned}\forall_{s\leq \frac{T}{n^{2\alpha}}} \sigma_Q B_{s}+\sigma_{A}n^{-\alpha}W^{(n,k)}_{sn^{2\alpha}}\in[a+\varepsilon, b-\varepsilon],
\\ \sigma_Q B_{\frac{T}{n^{2\alpha}}}+\sigma_{A}n^{-\alpha}W^{(n,k)}_{T}\in[a''+\varepsilon,b''-\varepsilon] \end{aligned} \Bigg|W\right)-e^{-n^\upsilon}\right].
\eeqnn
Similarly, for the last segment (from time $f'(n)-f_K(n)$ to time $f'(n)$), we also have, almost surely,
\beqnn
&&\inf\limits_{x\in [a'' n^\alpha,b'' n^\alpha]} \bfP_\mu\left(\begin{split}\forall_{i\leq f'(n)-f_K(n)} S_{f_K(n)+i}\in[an^\alpha, b n^\alpha],\\ S_{f'(n)}\in[a'n^\alpha,b'n^\alpha],\xi_{f_K(n)+i}\leq r_n \end{split} \Bigg|S_{f_K(n)}=x\right)
\\&\geq& 1_{A_n}\Big[\inf\limits_{x\in [a'',b'']}\bfP^x\Big(\forall_{s\leq \frac{T}{n^{2\alpha}}} \sigma_Q B_{s}+\sigma_{A}n^{-\alpha}W^{(n,k)}_{sn^{2\alpha}}\in[a'+\varepsilon, b'-\varepsilon]|W\Big)-e^{-n^\upsilon}\Big].
\eeqnn
Now it is time to introduce $G_{n}.$  Denote $G_{n}:=\cap_{k=0}^{K}J_{n,k}, $ where
$$J_{n,k}:=\{\frac{1}{2}e^{-X_{n,k}}\geq e^{-n^{\upsilon}}\}, \forall k<K,~~~ J_{n,K}:=\{\frac{1}{2}e^{-Y_{n,K}}\geq e^{-n^{\upsilon}}\},$$
\beqnn X_{n,k}&:=&-\ln\inf\limits_{x\in [a_0,b_0]} \bfP^x\left(\begin{aligned}\forall_{s\leq \frac{T}{n^{2\alpha}}} \sigma_Q B_{s}+\sigma_{A}n^{-\alpha}W^{(n,k)}_{sn^{2\alpha}}\in[a+\varepsilon, b-\varepsilon],
\\ \sigma_Q B_{\frac{T}{n^{2\alpha}}}+\sigma_{A}n^{-\alpha}W^{(n,k)}_{T}\in[a''+\varepsilon,b''-\varepsilon] \end{aligned} \Bigg|W\right),\eeqnn
$$Y_{n,K}:=-\ln\inf\limits_{x\in [a'',b'']} \bfP^x\Big(\forall_{s\leq D} \sigma_Q B_{s}+\sigma_{A}n^{-\alpha}W^{(n,K)}_{sn^{2\alpha}}\in[a'+\varepsilon, b'-\varepsilon]|W\Big).$$

Now we show $\lim\limits_{n\rightarrow+\infty}\mathbf{1}_{G_n}=1$~~a.s. by using Lemma 3.1.
Note that $$e^{-X_{n,k}}\geq e^{-X_{n,k,D-1}}e^{-X'_{n,K,1}},$$ where 
\beqnn e^{-X_{n,k,D-1}}&:=&\inf\limits_{x\in [a_0,b_0]}\bfP^x\left(\begin{aligned}\forall_{s\leq D-1}~ \sigma_Q B_{s}+\sigma_An^{-\alpha}W^{(n,k)}_{sn^{2\alpha}}\in[a+\varepsilon, b-\varepsilon],
\\ \sigma_Q B_{D-1}+\sigma_A n^{-\alpha}W^{(n,k)}_{_{(D-1)n^{2\alpha}}}\in[a''',b''']\end{aligned}\Bigg|W\right),\eeqnn
\beqnn e^{-X'_{n,k,1}}:=&\inf\limits_{x\in [a''',b''']}&\bfP^x\left(\begin{aligned}\forall_{0\leq s\leq 1}~ \sigma_A\left(W^{(n,k)}_{_{(D-1+s)n^{2\alpha}}}-
W^{(n,k)}_{_{(D-1)n^{2\alpha}}}\right)
\\+\sigma_QB_{s}\in[a''+\varepsilon, b''-\varepsilon]\end{aligned}\Bigg|W\right).\eeqnn
When $n$ is large enough, we can see $\{X_{n,k,D-1}\leq n^{\frac{\upsilon}{2}}\}\cap\{X'_{n,k,1}\leq n^{\frac{\upsilon}{2}}\}\subseteq J_{n,k}.$ Moreover, noting that for fixed $n$ and any $k<K,$ we have $X_{n,k,D-1}\stackrel{d}{=}X_{n,0,D-1}, X'_{n,k,1}\stackrel{d}{=}X'_{n,0,1} ,$ where ``$\stackrel{d}{=}$'' stands for equidistribution. Therefore, we have
\beqnn\bfP(|1_{G_{n}}-1|\geq \varepsilon)=\bfP(G^c_{n})&\leq& K\bfP(J^c_{n,0})+\bfP(J^c_{n,K})
\\&\leq &K[\bfP(X_{n,0,D-1}> n^{\frac{\upsilon}{2}})+\bfP(X'_{n,0,1}> n^{\frac{\upsilon}{2}})] + \bfP(J_{n,K}).
\eeqnn
By Lemma 3.1, we have  $\sum_{i=1}^{+\infty}\bfP(|1_{G_{n}}-1|\geq \varepsilon)\leq \sum_{i=1}^{+\infty}3Kn^{-q}<+\infty.$
According to the Borel-Cantelli 0-1 law, we know $\lim\limits_{n\rightarrow+\infty}\mathbf{1}_{G_{n}}=1.$ a.s., hence $\lim\limits_{n\rightarrow+\infty}\mathbf{1}_{A_{n}}=1,$ a.s..
Considering the definition of event $G_{n},$ for large enough $n$, we have
\beqnn
&&\inf\limits_{x\in [a_0 n^\alpha,b_0 n^\alpha]} \bfP_\mu\left(\forall_{f(n)\leq i\leq f(n)+n} \frac{S_{i}}{n^\alpha}\in[a, b],\frac{S_{f(n)+n}}{n^\alpha}\in[a',b'],\xi_i\leq r_n|S_{f_k(n)}=x\right)
\\&\geq&\mathbf{1}_{A_n}\prod\limits_{k=0}\limits^{K-1}\left\{\frac{1}{2}e^{-X_{n,k,D-1}-X'_{n,k,1}}\right\}\times\frac{1}{2}e^{-Y_{n,K}}.
\eeqnn
Now we take the limit as follows.
\beqlb\label{sec-5}
&&\liminf\limits_{n\rightarrow+\infty}\frac{\ln\inf\limits_{x\in[a_0 n^{\alpha}, b_0 n^{\alpha}]}\bfP_\mu\Big(\forall_{f(n)\leq i\leq f'(n)} \frac{S_{i}}{n^\alpha}\in[a, b],\frac{S_{f'(n)}}{n^\alpha}\in[a',b'],\xi_i\leq r_n|S_{f_k(n)}=x\Big)}{n^{1-2\alpha}}\nonumber
\\&\geq &\liminf\limits_{n\rightarrow+\infty}\frac{\sum_{k=0}^{K-1}(-\ln 2-X_{n,k,D-1}-X'_{n,k,1})-\ln2-Y_{n,K}}{n^{1-2\alpha}}\nonumber
\\&\geq & \frac{-\ln 2+\bfE(\tilde{X}''_{D-1})+\bfE(\tilde{X}''_{1})}{D},
\eeqlb
where
$$\tilde{X}''_{1}:=-\ln\inf\limits_{x\in[a'',b'']} \bfP^x\big(\forall_{0\leq s\leq 1} \sigma_Q B_{s}+\sigma_A W_s\in [a'+\varepsilon,b'-\varepsilon] |W),$$
\beqnn\tilde{X}''_{D-1}:=-\ln\inf\limits_{x\in [a_0,b_0]}\bfP^x\left(\begin{aligned}\forall_{s\leq D-1}~ \sigma_Q B_{s}+\sigma_A W_s\in[a+\varepsilon, b-\varepsilon],\\
 \sigma_Q B_{D-1}+\sigma_A W_{D-1}\in[a''',b''']\end{aligned}\Bigg|W\right).\eeqnn
The last inequality in (3.6) is obtained by using the Borel-Cantelli 0-1 law just like the proof of \cite[Corollary 2.1]{LY2018}. Moreover, By \cite[Theorem 3.1]{LY2018}, we know $\bfE(\tilde{X}''_{1})<+\infty.$
According to the  $L^1$ convergence in \cite[Theorem 2.1]{LY2018} and the continuity of $\gamma$, we get the final result of lower bound by taking $D\rightarrow +\infty,$ $\varepsilon\rightarrow 0$ (in this order).

\bigskip

{\bf Proof of Theorem 2.1: the upper bound}

The upper bound can be obtained by similar method as the lower bound. Here we let $T:=\lceil Dn^{2\alpha}\rceil, K:=\lfloor\frac{n}{T}\rfloor.$
Adjusting the definition of $G_n$ to $$G_n:=\bigcap_{k=0}^{K-1}\Big\{\sup_{x\in\bfR}\bfP^x\Big(\forall_{s\leq D} \sigma_Q B_{s}+\sigma_{A}n^{-\alpha}W^{(n,k)}_{sn^{2\alpha}}\in[a-\varepsilon, b+\varepsilon]\big|W\Big)\geq e^{-n^{\upsilon}}\Big\}.$$
Other marks and definitions will not need to be changed. It is true that
\beqnn
&&\ln\sup\limits_{x\in \bfR}\bfP_\mu
(\forall_{f(n)\leq i\leq f(n)+n} an^\alpha\leq S_{i}\leq bn^\alpha|S_{f(n)}=x)
\\&\leq&\mathbf{1}_{A_n}\ln\sup\limits_{x\in \bfR}\bfP_{\mu}
(\forall_{f(n)\leq i\leq f(n)+n} x+(U_i-U_{f(n)})+(M_i-M_{f(n)}) \in[an^\alpha,bn^\alpha])
\\&\leq&\mathbf{1}_{A_n}\ln\prod\limits_{k=0}\limits^{K-1}\sup\limits_{x\in\bfR} \bfP_\mu(\forall_{i\leq T} x+U^{(n,k)}_{i}+\sigma_{A}W^{(n,k)}_{i}\in[an^\alpha-2n^\rho, bn^\alpha+2n^\rho]|W).
\eeqnn
On the event $A_n,$ for any $k=0,1,...K-1$ we have
\beqnn
&&\bfP_\mu(\forall_{i\leq T} x+U^{(n,k)}_{i}+\sigma_{A}W^{(n,k)}_{i}\in[an^\alpha-2n^\rho, bn^\alpha+2n^\rho]|W)
\\&\leq& \bfP^{x}_{\mu}(\forall_{i\leq T}B_{\Gamma^{(n,k)}_{i}}+\sigma_{A}W^{(n,k)}_{i}\in[an^\alpha-3n^\rho, bn^\alpha+3n^\rho]|W)+\bfP_{\mu}(C^{c}_{n,k})
\\&\leq& \bfP^{x}(\forall_{s\leq T}B_{s\sigma^2_Q}+\sigma_{A}W^{(n,k)}_{s}\in[an^\alpha-5n^\rho, bn^\alpha+5n^\rho]|W)
+\bfP_{\mu}(C^c_{n,k})+\bfP(D^c_{n})
\\&\leq& \bfP^x\big(\forall_{s\leq T} \sigma_Q B_{s}+\sigma_{A}W^{(n,k)}_{s}\in[an^\alpha-5n^\rho, bn^\alpha+5n^\rho]|W\big)+e^{-n^{\upsilon}}.
\eeqnn
Noting that $\frac{T}{n^{2\alpha}}\geq D,$ by the scaling property of Brownian motion we obtain
\beqnn&&\ln\sup\limits_{x\in\bfR}\bfP^x_\mu
(\forall_{f(n)\leq i\leq f(n)+n} an^\alpha\leq S_{i}\leq bn^\alpha|S_{f(n)}=x)
\\&\leq&\mathbf{1}_{A_n}\left\{\sum\limits_{k=0}\limits^{K-1}\ln\Big[\sup\limits_{x\in\bfR} \bfP^x\Big(\forall_{s\leq D} \sigma_Q B_{s}+\sigma_{A}n^{-\alpha}W^{(n,k)}_{sn^{2\alpha}}\in[a-\varepsilon, b+\varepsilon]|W\Big)+e^{-n^{\upsilon}}\Big]\right\}
\\&\leq&\mathbf{1}_{A_n}\left\{\sum\limits_{k=0}\limits^{K-1}\ln\Big[2\sup\limits_{x\in\bfR} \bfP^x\big(\forall_{s\leq D} \sigma_Q B_{s}+\sigma_{A}n^{-\alpha}W^{(n,k)}_{sn^{2\alpha}}\in[a-\varepsilon, b+\varepsilon]|W\big)\Big]\right\}.
\eeqnn
Taking $n\rightarrow +\infty,$ just like the instruction of the corresponding part in the proof of the lower bound, we get
\beqnn&&\limsup_{n\rightarrow+\infty}\frac{\ln\sup\limits_{x\in\bfR}\bfP^x_\mu
\left(\forall_{f(n)\leq i\leq f(n)+n} an^\alpha\leq S_{i}\leq bn^\alpha|S_{f(n)}=x\right)}{n^{1-2\alpha}}
\\&\leq&\frac{\ln2+\bfE\big[ \bfP^x_{\mu}\big(\forall_{s\leq D} \sigma_Q B_{s}+\sigma_{A}W_s\in[a-\varepsilon, b+\varepsilon]|W\big)\big]}{D}.
\eeqnn
At last, let $D\rightarrow+\infty,$ $\varepsilon\rightarrow 0.$  Applying the $L^1$ convergence in \cite[Theorem 2.1]{LY2018}, we complete the proof.

{\bf Proof of Corollary 2.1}

Now we turn to prove the corollary. The proof is mainly depending on Markov property and the continuity of the functions $g$ and $h$.
Here we only write down the proof of the lower bound. The method for the upper bound is as the same way and easier so we left it to readers.
Choose $m\in\bfN$ and let $d_n=\lfloor\frac{n}{m}\rfloor.$ First we try to change the length of time $n$ to $d_n$ in Theorem 2.1. Without loss of generality,
assume that $a<a_0<a'<0<b'<b_0<b.$ Choose $v>0$ such that $(\frac{m+v}{m})^{\alpha}<\min\{\frac{a}{a_0},\frac{b}{b_0}\}.$ Let $n$ large enough to satisfy
that $n\leq (m+v)d_n,$ then
\beqnn
&&\inf\limits_{x\in[a_0 n^{\alpha}, b_0 n^{\alpha}]}\bfP_\mu
\left(\forall_{f(n)\leq i\leq f(n)+d_n} \frac{S_{i}}{n^\alpha}\in[a, b], \frac{S_{f(n)+d_n}}{n^\alpha}\in[a', b'],\xi_i\leq r_n\Big|S_{f(n)}=x\right)
\\&\geq&\inf\limits_{x\in[a_0(m+v)^{\alpha} N^{\alpha}, b_0(m+v)^{\alpha} N^{\alpha}]}\bfP_\mu
\left(\begin{aligned}\forall_{f(n)\leq i\leq f(n)+d_n} N^{-\alpha}S_{i}\in[am^\alpha, bm^\alpha],\\
N^{-\alpha}S_{f(n)+d_n}\in[a'm^\alpha, b'm^\alpha],~\xi_i\leq r_n\end{aligned}\Bigg|S_{f(n)}=x\right).
\eeqnn
Therefore, we have \beqlb\label{sect-7}&&\liminf_{n\rightarrow +\infty}\frac{\ln\inf\limits_{x\in[a_0 n^{\alpha}, b_0 n^{\alpha}]}\bfP_\mu
\left(\begin{aligned}\forall_{f(n)\leq i\leq f(n)+d_n} S_{i}\in[a n^\alpha, bn^\alpha], \\ S_{f(n)+d_n}\in[a'n^\alpha, b'n^\alpha],\xi_i\leq r_n\end{aligned}\Bigg|S_{f(n)}=x\right)}{n^{1-2\alpha}}\nonumber
\\&\geq&-m^{2\alpha-1}\frac{\sigma^2_Q}{(bm^{\alpha}-am^{\alpha})^2}\gamma(\frac{\sigma_A}{\sigma_Q})
=-\frac{1}{m}\frac{\sigma^2_Q}{(b-a)^2}\gamma(\frac{\sigma_A}{\sigma_Q}).
\eeqlb
Let $\bar{\varepsilon}=\min\left\{\frac{\inf_{s\in[0,1]}(h(s)-g(s))}{6}, \frac{b_0-a_0}{2},\frac{b'-a'}{2}, h(0)-b_0, a_0-g(0),h(1)-b', a'-g(1)\right\}.$ Choose an $\varepsilon\in (0,\bar{\varepsilon}),$ then there exists a function $l$ defined on $[0,1]$ such that $l(0)\in(a_0+\varepsilon,b_0-\varepsilon), l(1)\in(a'+\varepsilon,b'-\varepsilon), l(t)\in(g(t)+3\varepsilon,h(t)-3\varepsilon)$ for every $t\in(0,1).$ Now choose $m\in\bfN$ large enough to satisfy that $\sup\limits_{0\leq|t-s|\leq2m^{-1}}\{|g(t)-g(s)|+|h(t)-h(s)|\}\leq \varepsilon. $
 Notice that adding the event $\xi_i\leq r_n$ will not bring any extra difficulty since $\{\xi_i\}_{i\in\bfN}$ is an independent random sequence in space $\bfP_\mu$.  Hence in the next inequalities, we omit the event $\xi_i\leq r_n$ for the simplicity of writing.
Denote $J_{k,n}=\big[(l(k/m)-\varepsilon)n^{\alpha},(l(k/m)+\varepsilon)n^{\alpha}\big],$
$f(n,k):=f(n)+kd_n.$ We have
\beqnn
&&\inf\limits_{x\in [a_0 n^\alpha, b_0 n^\alpha]}\bfP_{\mu}\left(\begin{aligned}\forall_{0\leq j\leq n} ~g\big(j/n\big)n^\alpha\leq S_{f(n)+j}\leq h\big(j/n\big)n^\alpha,\\ a'n^\alpha\leq S_{f(n)+n}\leq b'n^\alpha~~~~~ \end{aligned}\Bigg|S_{f(n)}=x\right)
\\&\geq& \inf\limits_{x\in [a_0 n^\alpha, b_0 n^\alpha]}\bfP_{\mu}\left(\forall_{0\leq j\leq d_{n}} g\Big(\frac{j}{n}\Big)n^\alpha\leq S_{f(n)+j}\leq h\Big(\frac{j}{n}\Big)n^\alpha, S_{d_{n}}\in J_{1,n}\Big|S_{f(n)}=x\right)
\\&\times&\prod\limits_{k=1}\limits^{m-1}\inf\limits_{x\in J_{k,n}}\bfP_{\mu}\left(\forall_{0\leq j\leq d_{n}}\begin{aligned} I_{j,n} , S_{f(n,k+1)}\in J_{k+1,n}\end{aligned}|S_{f(n,k)}=x\right)
\\&\times&\inf\limits_{x\in J_{m,n}}\bfP_{\mu}(\forall_{0\leq j\leq d_{n}} a'n^\alpha\leq S_{f(n,m)+j}\leq b'n^\alpha|S_{f(n,m)}=x\Big),
\eeqnn where the event $$I_{j,n,k}:=\left\{g\Big(\frac{k d_n+j}{n}\Big)\leq \frac {S_{f(n)+kd_n+j}}{n^\alpha}\leq h\Big(\frac{k d_n+j}{n}\Big)\right\}.$$
Denote $$\underline{h}_{a,m}:=\inf_{x\in\left[0\vee\frac{a-1}{m},1\wedge\frac{a+1}{m}\right]}h(x),~~~\overline{g}_{a,m}:=\sup_{x\in\left[0\vee\frac{a-1}{m},1\wedge\frac{a+1}{m}\right]}g(x).$$
Then we have
\beqnn
&&\liminf_{n\rightarrow +\infty}\frac{\ln\inf\limits_{x\in [a_0 n^\alpha, b_0 n^\alpha]}\bfP_{\mu}\left(\forall_{0\leq j\leq n} I_{j,n,0},~S_{f(n)+n}\in[a'n^\alpha,b'n^\alpha]|S_{f(n)}=x\right)}{n^{1-2\alpha}}
\\&\geq& \frac{1}{m}\left[\sum_{a=0}^{m-1}\frac{\sigma^2_Q}{(\underline{h}_{a,m}-\overline{g}_{a,m})^2}\gamma\left(\frac{\sigma_A}{\sigma_Q}\right)
+\frac{\sigma^2_Q}{(b'-a')^2}\gamma\left(\frac{\sigma_A}{\sigma_Q}\right)\right].
\eeqnn
Hence we complete the proof by taking $m\rightarrow +\infty.$ 

\ack
I want to thank my supervisor Wenming Hong for his constant concern on my
work and giving me a good learning environment. I also want to thank Bastien Mallein for giving me a lot of valuable advice and useful tips, especially for his help in \cite{LY2018}, which is the basis and preparation work for this paper.


\begin{thebibliography}{99}

\bibitem{GHS2011}
N.Gantert, Y.Hu, Z.Shi. Asymptotics for the survival probability in a killed branching random walk. Ann. Probab. Statist. (47) 111-129. 2011.
\bibitem{Lif2000}
M. Lifshits. Lecture notes on strong approximation. 2000.
\bibitem{L2006}
M. A. Lifshits.~ Bibliography of small deviation probabilities (available via www.proba.jussieu.fr/pageperso/smalldev/biblio.html), 2006.
\bibitem{LY2018}
Y. Lv. Brownian motion between two random trajectories. Accepted to Markov Processes and Related Fields. ArXiv e-prints, arXiv:1802.03876v4, 2018.
\bibitem{LY201809}
Y. Lv. Branching random walk in random environment with random absorption wall. ArXiv e-prints, arXiv:1809.04969, 2018.
\bibitem{Mal2015}
B. Mallein.~ Maximal displacement of a branching random walk in time-inhomogeneous environment. Stochastic Process. Appl. (125) 3958-4019. 2015.
\bibitem{MM2015}
B. Mallein, P. Mi{\l}o\'{s}.~ Brownian motion and random walks above quenched random wall.  Accepted to Ann. Inst. Henri Poincar\'{e} Probab. Stat. ArXiv e-prints, arXiv:1507.08578, 2015.
\bibitem{MM2016}
B. Mallein, P. Mi{\l}o\'{s}.~ Maximal displacement of a supercritical branching random walk in a time-inhomogeneous random environment.  Accepted to Stochastic Process. Appl. ArXiv e-prints, arXiv:1507.08835,  2016.
\bibitem{Mog1974}
A.A.Mogul'ski\v{\i}.  Small deviations in the space of trajectories. Theory Probab. Appl. 19, 726-736, 1974.
\bibitem{QM1995}
Q. Shao. A Small Deviation Theorem for Independent Random Variables. Theory Probab. Appl. 40(1):225¨C235, 1995.
\end{thebibliography}
\end{document}